\documentclass[12pt]{amsart}
\usepackage{amssymb}

\theoremstyle{plain}
\newtheorem*{thm}{Theorem}
\newtheorem*{lem}{Lemma}
\newtheorem*{prop}{Proposition}

\theoremstyle{remark}
\newtheorem*{rem}{Remark}
\newtheorem*{ack}{Acknowledgment}
\newtheorem*{notat}{Notations and Conventions}

\def\hgt{\operatorname{height}}
\def\ann{\operatorname{ann}}
\def\nio{\operatorname{nio}}
\def\Z{\mathfrak{Z}}
\def\clKdim{\operatorname{clKdim}}
\def\kdim{\operatorname{Kdim}}
\def\*{\!*\!}
\def\G{\Gamma}

\begin{document}

%Topmatter%%%%%%%%%%%%%%%%%%%%%%%%%%%%%%%%

\title{Polycyclic-by-finite group algebras are catenary}

\author{Edward~S. Letzter}
\address{Department of Mathematics\\
 Texas A\&M University\\
 College Station, TX 77843}
\email{letzter@math.tamu.edu}
\thanks{Research of the first author supported in part by NSF Grant
DMS-9623579}

\author{Martin Lorenz}
\address{Department of Mathematics\\
        Temple University\\
        Philadelphia, PA 19122-6094}
\email{lorenz@math.temple.edu}
\thanks{Research of the second author supported in part by NSF Grant
DMS-9618521}



\begin{abstract}
We show that group algebras $kG$ of polycyclic-by-finite groups $G$,
where $k$ is a field, are catenary: If $P = I_0 \subsetneq I_1
\subsetneq \cdots \subsetneq I_m = P'$ and $P = J_0 \subsetneq J_2
\subsetneq \cdots \subsetneq J_n = P'$ are both saturated chains of
prime ideals of $kG$, then $m = n$.
\end{abstract}

\maketitle

%end Topmatter%%%%%%%%%%%%%%%%%%%%%%%%%%%%%%%%%%

\section*{Introduction} \label{intro}

Group algebras $kG$ of polycyclic-by-finite groups $G$, over a field
$k$, are among the most prominent special classes of noncommutative
noetherian rings, partly because the only known noetherian group
algebras (over fields) are of the form $kG$, and partly because of the
intrinsic group theoretic interest in $G$. Indeed, the investigation
of the prime ideals of $kG$ was a major enterprise throughout the
1970's. Zalesskii's work \cite{Z} initiated a program that was
essentially completed in 1980 by \cite{LP}, \cite{L} and \cite{B}; the
most notable contribution was by Roseblade \cite{R}. The
aforementioned papers resulted in a detailed, albeit complicated,
description of the prime ideals of $kG$, and the subject has been
largely dormant ever since. Disappointingly, the question of whether
$kG$ is, or is not, catenary (see below) remained unanswered, despite
some attempts in that direction: \cite[Corrigenda]{R} gave an
affirmative answer for $\Z$-groups and \cite{L} did the same for
orbitally sound-by-(finite nilpotent) groups. But the intervening
almost 20 years have witnessed substantial progress in our
understanding of prime ideals in both noetherian rings and their
extensions (see, e.g., \cite{GW}), and so we can now offer in the
present paper a proof of the following

\begin{thm}
Group algebras of polycyclic-by-finite groups are catenary.
\end{thm}

Recall that a ring $R$ is called catenary when the following
holds: For every inclusion $P \subsetneq P'$ of prime ideals in $R$,
if $P = I_0 \subsetneq I_1 \subsetneq \cdots \subsetneq I_m = P'$ and
$P = J_0 \subsetneq J_2 \subsetneq \cdots \subsetneq J_n = P'$ are
both chains of prime ideals that cannot be lengthened by the insertion
of extra primes, then $m = n$.  A fundamental result of commutative
algebra asserts that all affine commutative $k$-algebras are catenary.
Various noncommutative generalizations of this result are known:
Affine PI-algebras \cite{S}, enveloping algebras of finite dimensional
solvable Lie (super-)algebras \cite{Ga}, \cite{Le} and, most recently,
certain algebras arising from quantum groups \cite{GL} have all been
proven to be catenary.

%===================================================

\begin{notat}
Throughout this note, $k$ will denote a commutative
field. Furthermore, $\subsetneq$ denotes strict inclusion, while
equality is permitted in $\subset$; similarly for $\supsetneq$ and
$\supset$.  We will write $\ann X_R$ when referring to the right
annihilator of $X$ in $R$ and $\ann _SY$ when referring to the left
annihilator of $Y$ in $S$.
\end{notat}

%%%%%%%%%%%%%%%%%%%%%%%%%%%%%%%%%

\section{Preliminaries} \label{prelim}

%=================================
\subsection{Catenarity} \label{cat}
A ring $R$ is called \textit{catenary} if, given prime ideals $Q
\subsetneq Q'$ of $R$, any two saturated chains of primes between $Q$
and
$Q'$ have the same length. Here, a chain $Q_0 \subsetneq Q_1 \subsetneq
\cdots \subsetneq Q_m$ of prime ideals of $R$ is \textit{saturated} if
$Q_{i-1}$ and $Q_i$ are adjacent for all $1 \leq i \leq m$; a pair of
prime ideals $Q$ and $Q'$ in $R$ are \textit{adjacent} if $Q \subsetneq
Q'$ and if there exist no prime ideals $Q''$ of $R$ such that $Q
\subsetneq Q'' \subsetneq Q'$.

Recall that the \textit{height} of a prime ideal $P$ of $R$ is the
supremum of the lengths $n$ of all chains $P=P_0\supsetneq
P_1\supsetneq\dots \supsetneq P_n$ of prime ideals descending from $P$.
In
all known noetherian rings the prime ideals have finite height; for
group algebras of polycyclic-by-finite groups, see \cite[Theorem
H1]{R}.  Assuming finiteness of heights, $R$ is certainly catenary if
the following condition is satisfied:
\begin{equation} \label{E:ht}
\hgt P' =\hgt P + 1 \text{ holds whenever $P\subsetneq P'$ are
adjacent prime ideals.}
\end{equation}
If every prime of $R$ contains \textit{exactly one} minimal prime,
then \eqref{E:ht} is actually equivalent with catenarity.

For group algebras $kG$ of polycyclic-by-finite groups $G$, the latter
condition on minimal primes holds. Indeed, if $F$ denotes the
finite radical of $G$, that is, the largest finite normal subgroup of
$G$, then the unique minimal prime of $kG$ contained in a given prime
ideal $P$ of $kG$ is $(P\cap kF)kG$; see \cite[Corollary 22]{R}. Thus,
our objective will be to verify \eqref{E:ht} for $kG$.

%=============================================
\subsection{Heights} \label{heights}
Let $kG$ be the group algebra of a polycyclic-by-finite group $G$.
Following Roseblade \cite{R} we put, for any ideal $I$ of $kG$,
$$
I^\dag = G\cap (1+I)\ .
$$
Thus $I^\dag$ is the kernel of the natural homomorphism of $G$
into the group of units of $kG/I$, and $I$ contains the kernel
of the natural homomorphism $kG\to k[G/I^\dag]$. The ideal $I$
is said to be \textit{faithful} iff $I^\dag = 1$. The image of
$I$ in $k[G/I^\dag]$ is always faithful.

The following Lemma will be used to reduce the proof of \eqref{E:ht}
to the case where $P$ is faithful. Its essence is in \cite[2.4]{R},
at least for Roseblade's $\Z$-groups; the general case follows
rather routinely. The term
$$
p_G(N)
$$
in the formula below denotes the \textit{$G$-plinth length} of the
normal subgroup $N$ of $G$. For the definition and basic properties
of plinth lengths, the reader is referred to \cite[2.3]{R}.
We remark that \cite[2.3]{R} assumes $N$ to be
polycyclic, but the notion of plinth length is insensitive to finite
factors and
applies to polycyclic-by-finite groups as well.

\begin{lem} Let $kG$ be the group algebra of a polycyclic-by-finite
group $G$,
let $P$ be a prime ideal of $kG$, and let $N$ be a normal subgroup of
$G$
that is contained in $P^\dag$. Then
$$
\hgt P = \hgt\overline{P} + p_G(N)\ ,
$$
where $\overline{P}$ denotes the image of $P$ in $k[G/N]$.
\end{lem}

\begin{proof}
Let $G_0$ denote a normal $\Z$-group of finite index in $G$ (see
\cite[1.3]{R})
and put $N_0=N\cap G_0$. Further, write
\begin{equation} \label{cutdown}
P\cap kG_0 = \bigcap_{g\in G} P_0^g
\end{equation}
for some prime ideal $P_0$ of $kG_0$; cf.~\cite[4.1]{R} or
\cite[Lemmas 1.3(i), 1.4(i)]{LP}.  Then $P^\dag\cap G_0\subset
P_0^\dag$, and so $N_0\subset P_0^\dag$. Letting
$\overline{\phantom{x}}$ denote images in $k[G/N]$, it follows that
$\overline{P}\cap k[\overline{G_0}] = \bigcap_{x\in\overline{G}}
\overline{P_0}^x$ and $\overline{P_0}$ is prime in
$k[\overline{G_0}]$. Thus, by \cite[8.1]{R} or \cite[Lemma
2.5(i)]{LP}, $\hgt P = \hgt P_0$ and $\hgt\overline{P} =
\hgt\overline{P_0}$.  Finally, $p_G(N)=p_{G_0}(N_0)$, since $G_0$ has
finite index in $G$.  Therefore, replacing $G$ by $G_0$ etc., we may
assume that $G$ is a $\Z$-group. In this case, the desired formula
follows from Roseblade's height formula in \cite[2.4]{R}, noting that
Roseblade's $\lambda$ equals $\hgt$ for $\Z$-groups and that
$p_G(\,.\,)$ is additive on short exact sequences of $G$-groups.
\end{proof}

%=============================================
\subsection{Second layer theory} \label{slc}
We review some
relevant results and background information from noetherian ring
theory.\smallskip

(i) Definitions of the right and left strong second layer conditions
of Jategaonkar may be found in \cite[p.~183]{GW}; also see \cite{J}.
When we
refer to noetherian rings having the strong second layer condition,
abbreviated \textit{sslc}, we mean rings satisfying both the right and
left strong second layer conditions.\smallskip

(ii) Group rings of polycyclic-by-finite groups over commutative
noetherian coefficient rings satisfy the sslc; see
\cite[A.4.6]{J}.\smallskip

(iii) Let $U$ and $V$ be noetherian rings satisfying the sslc. Suppose
that $B$ is a $U$-$V$-bimodule finitely generated and faithful on each
side. It was proved by Jategaonkar that $U$ and $V$ have the same
classical Krull dimension; see, for example, \cite[12.5]{GW} or
\cite[8.2.8]{J}.
We will use $\clKdim (\, . \,)$ to denote classical Krull
dimension.\smallskip

(iv) Let $U$ be a prime noetherian ring, and let $M$ be a nonzero left
$U$-module. Recall that $u \in U$ is \textit{regular} provided $uv$ and
$vu$ are nonzero for all nonzero $v \in U$. Further recall that $M$ is
\textit{torsion free} provided $u.m \ne 0$ for every regular element $u
\in U$ and every $0\ne m \in M$. It follows from Goldie's Theorem,
when $_UM$ is torsion free, that every nonzero submodule of $M$ is
faithful. If $N$ is a $U$-module for which every nonzero submodule is
faithful, then we say that $N$ is {\it fully faithful} over
$U$.\smallskip

(v) Let $P$ and $Q$ be prime ideals in a noetherian ring $U$
satisfying the sslc. We say that $P$ is \textit{linked} to $Q$, and
write $P \rightsquigarrow Q$, if there exists a $U$-$U$-bimodule
factor $B$ of $(P\cap Q)/PQ$ with the following properties: $\ann_UB =
P$, $\ann B_U = Q$, and $B$ is torsion free on each side as a
$U/P$-$U/Q$-bimodule. A graph structure is thus imposed on the set of
prime ideals, and if $X$ is a subset thereof then the \textit{link
closure} of $X$ will be defined to be the union of the connected
components of the prime ideals in $X$. Note, by (iii), that if $P
\rightsquigarrow Q$ then $\clKdim (U/P) = \clKdim (U/Q)$.\smallskip

(vi) Continue to let $U$ be a noetherian ring satisfying the sslc, and
let $M$ be a finitely generated left $U$-module. A prime ideal $P$ of
$U$ is said to be an \textit{associated prime} of $M$ if there exists
a submodule $N$ of $M$ such that $\ann _UN = P$ and such that $N$ is a
fully faithful $U/P$-module. By noetherianity, the set of associated
primes of $M$ is not empty. Also, if $M'$ is an essential extension of
$M$, then the sets of associated primes of $M$ and $M'$ are the same.
Let $X$ denote the link closure of the set of associated primes of
$M$. In \cite[11.4]{GW}, for example, it is proved that $M$ is
annihilated
by a (finite) product of prime ideals in $X$. (This last assertion
follows from an iterated application of the Main Lemma of Jategaonkar;
see \cite{J}.)

%=============================================
\subsection{Krull dimension} \label{krull}
Here we collect some facts on \textit{Krull dimension}, in the sense
of Gabriel and Rentschler, that will be needed later. The Krull
dimension of a module will be denoted $\kdim(M)$; for the definition
and basic facts concerning $\kdim$, we refer the reader to
\cite[Chapter 13]{GW}. All rings and modules discussed below are
assumed noetherian, and so $\kdim(M)$ will always exist. Further,
$\kdim(U)$, for a ring $U$, will stand for $\kdim(U_U)$.
\smallskip

(i) A nonzero module $M$ is called \textit{homogeneous} (or
\textit{$\alpha$-homogeneous}, if $\kdim(M)= \alpha$) provided
all nonzero submodules of $M$ have the same Krull dimension as $M$.
Prime rings are homogeneous; see \cite[p.~229]{GW}. Further, finite
subdirect products and extensions of $\alpha$-homogeneous modules are
easily seen to
be $\alpha$-homogeneous.\smallskip

(ii) Assume the ring $U$ is (right) homogeneous. If the element $u\in
U$ satisfies $\kdim(U/uU)<\kdim(U)$ then $u$ is regular. This
assertion can be proved as follows: First, the left regularity of $u$
follows from \cite[6.8.12]{McCR}. As to right regularity (which is due
to Stafford \cite{St}), observe that $(U/uU)_U$ maps onto
$(u^{n-1}U/u^nU)_U$
for all $n>0$. Hence, $\kdim(u^{n-1}U/u^nU)_U\le\kdim(U/uU)_U<\kdim(U)$
and, consequently, $\kdim(U/u^nU)_U<\kdim(U)$. But, by ``Fitting's
Lemma"
(\cite[2.3.2(ii)]{McCR}), there is an $n>0$ such that the
right annihilator, $\ann(u^n)_U$, of $u^n$ in $U$ embeds into $U/u^nU$.
Therefore, $\kdim(\ann(u^n)_U)<\kdim(U)$, whence $\ann(u^n)_U=0$,
by homogeneity. This proves that $u$ is right regular.
\smallskip

(iii) For group algebras $kG$ of polycyclic-by-finite groups $G$, one
has $\kdim(kG)=h(G)$, the Hirsch number of $G$; see \cite{Sm}.
If $H$ is a subgroup of finite index in $G$, then for any (right)
$kG$-module $M$, one has
$$
\kdim(M_{kG})=\kdim(M_{kH})\ .
$$
If $H$ is normal, this is covered by a more general result on
finite normalizing extensions; see \cite[Th\'eor\`eme 5.3]{Lm}.
The general case follows by passing to the \textit{normal core}
$H_G=\bigcap_{g\in G}H^g$, a normal subgroup of finite index in $G$
(and in $H$).
\smallskip

(iv) Let $G$, $H$, and $M$ be as above and suppose that $M=WkG$ for some
submodule $W$ of $M_{kH}$. Then
$$
\kdim(M)=\kdim(W)\ .
$$
Indeed, by (iii), it suffices to prove $\kdim(M_{kN})=\kdim(W_{kN})$,
where $N=H_G$ is the normal core of $H$. But, as $kN$-module, $M=WkG$
is a finite sum of $G$-conjugates of $W_{kN}$ all of which have
isomorphic submodule lattices, and hence the same $\kdim$.  Thus,
$\kdim(M_{kN})=\kdim(W_{kN})$, as required.
\smallskip

(v) Keeping the above notations $G$ and $H$, a $kG$-module $M$ is
$\alpha$-homogeneous
if and only if it is so as $kH$-module.  To verify this
claim, first note that $\alpha$-homogeneity as $kH$-module certainly
entails $\alpha$-homogeneity over $kG$, by (iii). Conversely, if
$M_{kG}$ is $\alpha$-homogeneous and $W$ is a nonzero submodule of
$M_{kH}$, then $WkG$ is a nonzero $kG$-submodule of $M$, and so
$\alpha=\kdim(WkG)=\kdim(W)$, by (iv). Finally, if $M=WkG$ for
some submodule $W$ of $M_{kH}$, then $M$ is $\alpha$-homogeneous if
and only if $W$ is. Indeed, if $M=WkG$ is $\alpha$-homogeneous and
$X$ is a nonzero submodule of $W$, then $\alpha=\kdim(XkG)=\kdim(X)$,
by (iv). Conversely, assume that $W$ is $\alpha$-homogeneous. It
suffices to show that $M=WkG$ is $\alpha$-homogeneous as $kN$-module,
where $N=H_G$. But $M_{kN}$ is a finite sum of
$G$-conjugates of $W_{kN}$, and all these are $\alpha$-homogeneous.
Hence
$M_{kN}$ is $\alpha$-homogeneous.

%%%%%%%%%%%%%%%%%%%%%%%%%%%%%%%%%%%%

\section{Induced Ideals} \label{ind}

\textit{Throughout this section, $R\subset S$ is an extension of
rings.
Additional assumptions will be required later.}

\subsection{Definition and basic properties} \label{def}
For any ideal $I$ of $R$, we put
$$
I^S=\ann(S/IS)_S=\{s\in S\mid Ss\subset IS\}\ .
$$
In other words, $I^S$ is the largest ideal of $S$ that is contained in
$IS$. We refer to $I^S$ as the \textit{induced ideal of $I$}, and we
will say that ideals in $S$ of this general form are \textit{induced}
from $R$. In the case of group algebras $R=kH\subset S=kG$, for some
subgroup $H$ of $G$, one also writes $I^G$ in place of $I^S$ and calls
these ideals \textit{induced from $H$}.

The induction operator
$$
(\,.\,)^S: \text{Ideals of $R$}\to\text{Ideals of $S$}
$$
obviously preserves inclusions. Furthermore,
\begin{equation} \label{prod}
I_1^SI_2^S\subset (I_1I_2)^S
\end{equation}
holds for any two ideals $I_1, I_2$ of $R$. Indeed,
$I_1^SI_2^S\subset I_1SI_2^S\subset I_1I_2S$.

\subsection{Minimal covering primes} \label{min}
Equation \eqref{prod} has the following consequence.

\begin{lem}
Let $P$ be a prime ideal of $S$, and let $Q$ be an ideal of $R$
maximal among ideals whose induced ideals are contained in $P$.  Then
$Q$ is prime.
\end{lem}

\begin{proof}
Let $I_1$ and $I_2$ be ideals of $R$ properly containing $Q$. Then
$I_1^S$ and $I_2^S$ are not contained in $P$, and hence neither is
their product.  By \eqref{prod}, $(I_1I_2)^S$ is not contained in $P$,
whence $I_1I_2$ is not contained in $Q$.
\end{proof}

In the situation of the above Lemma, $P$ will often actually be minimal
over $Q^S$.
A sufficient condition for this to happen is given in the following

\begin{prop}
Suppose that $S$ is finitely generated and free as a left
$R$-module. Assume further that $R$ and $S$ are noetherian and satisfy
the sslc. Let $Q$ be a prime ideal of $R$, and let $P$ be a prime
ideal of $S$ containing $Q^S$. Then the following are equivalent:
\begin{enumerate}
\item[i.] $P$ is minimal over $Q^S$;
\item[ii.] $S/P$ and $R/Q$ have the same classical Krull dimension;
\item[iii.] $Q$ is maximal among ideals in $R$ whose induced ideal is
contained in $P$.
\end{enumerate}
\end{prop}

\begin{proof}
Set $M = S/QS$ and $J = Q^S = \ann M_S$. By hypothesis on $_RS$, $M$
is free of finite rank as left $R/Q$-module. In particular, it follows
from (\ref{slc}\,iii) that $\clKdim(R/Q) = \clKdim(S/J)$, and so
$\clKdim(S/P) \leq \clKdim(R/Q)$.\smallskip

(i $\Rightarrow$ ii) Note that $M$ is nonzero
and torsion free as a left $R/Q$-module. Therefore, by \cite[7.7]{GW},
there exist prime ideals $P_1,\ldots,P_m$ of $S$, and a series of
$R$-$S$-bimodules
$$
0 = M_0 \subsetneq M_1 \subsetneq \cdots \subsetneq M_m = M ,
$$
with the following properties: For each $1 \leq i \leq n$, the right
annihilator in $S$ of $M_i/M_{i-1}$ is equal to $P_i$, and
$M_i/M_{i-1}$ is faithful on each side as an
$R/Q$-$S/P_i$-bimodule. Note that all of the prime ideals in $S$
minimal over $J$, and in particular $P$, are included among
$P_1,\ldots,P_m$. Therefore, $\clKdim(S/P) = \clKdim(R/Q)$, by
(\ref{slc}\,iii).\smallskip

(ii $\Rightarrow$ i) If $P$ is not minimal over $J$, then
$$
\clKdim(S/P) < \clKdim(S/J) = \clKdim(R/Q) .
$$

(ii $\Rightarrow$ iii) Suppose that $I$ is an ideal of $R$, properly
containing $Q$, such that $I^S \subset P$. Since $S$ is free as
a left $R$-module, $\ann _R(S/IS) = I$. Therefore, by (\ref{slc}\,iii),
$$
\clKdim(R/Q) > \clKdim(R/I) = \clKdim(S/I^S) \geq \clKdim(S/P) .
$$

(iii $\Rightarrow$ ii) Suppose that $\clKdim(S/P) \ne \clKdim(R/Q)$,
and so $\clKdim(S/P) < \clKdim(R/Q)$. We will show that there exists
an ideal $I$ of $R$, strictly containing $Q$, such that $I^S \subset
P$.

By \cite[12.3]{GW}, there exists a series of $R$-$S$-bimodules
$$
0 \subset L_1 \subsetneq L_2 \subset M
$$
such that $\ann (L_2/L_1)_S = P$, and such that $L_2/L_1$ is torsion
free as a right $S/P$-module. Let $K$ be an $R$-$S$-sub-bimodule of
$M$, containing $L_1$, and maximal with respect to $K \cap L_2
\subset L_1$. Set $N = M/K$. Note that $L_2/L_1$ embeds as an
$R$-$S$-bimodule into $N$; we will let $L$ denote the image of this
embedding. Observe that every nonzero $R$-$S$-sub-bimodule of $N$
intersects $L$ nontrivially. Also, by (\ref{slc}\,iv), $L$ is a fully
faithful right $S/P$-module.

Let $Q'$ be an associated prime of $_RL$, and let $L' = \ann Q'
_L$. Note that $L'$ is a nonzero $R$-$S$-sub-bimodule of $L$, and that
$Q' = \ann_RL'$. In particular, $L'$ is an $R/Q'$-$S/P$-bimodule
finitely generated and faithful on each side. Consequently,
by (\ref{slc}\,iii),
$\clKdim(R/Q') = \clKdim (S/P)$. Now let $X$ denote the link closure of
the
associated primes of $_RL$. By (\ref{slc}\,v), if $Q'' \in X$, then
$\clKdim(R/Q'') = \clKdim (S/P) < \clKdim (R/Q)$.

Now let $T$ be a left $R$-submodule of $N$ maximal with respect to $T
\cap L = 0$. Set $V = N/T$ and $I = \ann _RV$. Note that $I \supset
Q$ and that $L$ embeds into $_RV$ as an essential submodule. It
now follows from (\ref{slc}\,vi) that $V$ is annihilated on the left by
a product $Q_1\cdots Q_n$ of prime ideals $Q_1,\ldots,Q_n \in
X$. However, the classical Krull codimension of each $Q_i$ is strictly
smaller than $\clKdim (R/Q)$, and so $I$ properly contains $Q$.

Since $IV = 0$, we have $IN \subset T$, and so $IN \cap L =
0$. However, $IN$ is an $R$-$S$-sub-bimodule of $N$, and so $IN = 0$
because each nonzero $R$-$S$-sub-bimodule of $N$ intersects $L$
nontrivially. Therefore, $IM \subset K$, and so $IM \cap L_2
\subset L_1$. Hence, $L_2/L_1$ is an
$R$-$S$-bimodule subfactor of $M/IM$. Since $\ann(L_2/L_1)_S=P$
and $M/IM\cong S/IS$, we conclude that
$$
P \supset \ann (M/IM)_S = \ann(S/IS)_S = I^S.
$$
The proposition follows.
\end{proof}

\begin{rem}
The argument establishing (iii $\Rightarrow$ ii) above is similar to
the proof in \cite[4.6]{Lt} of ``lying over'' for finite extensions of
noetherian rings with the second layer condition.
\end{rem}

%===============================
\subsection{The case of group rings} \label{grouprings}
We will broaden the list of equivalent properties in Proposition
\ref{min} by adding another condition concerning heights. This
strengthening of the proposition seems, however, to require further
restrictions on the ring embedding $R\subset S$.  Therefore, we will
concentrate here on extensions of group algebras,
$$
R=kH\subset S=kG\ ,
$$
where $G$ is a polycyclic-by-finite group and $H$ is a subgroup of
finite index in $G$.

We start with a lemma on (Krull) homogeneity whose essence is gleaned
from \cite[Lemma 4.2]{B}.

\begin{lem}
Let $G$ be a polycyclic-by-finite group, let $H$ be a subgroup of
finite index in $G$, and put $N=H_G$, the normal core of $H$. Further,
let $Q$ be a prime ideal of $kH$ and put $I=Q^G\cap kN$. Then
$$
kH/Q,\ kG/Q^G, \text{ and }\ kG/IkG
$$
are all (right) $\alpha$-homogeneous rings, where $\alpha=\kdim(kH/Q)$.
\end{lem}

\begin{proof}
First, $kH/Q$ is prime and therefore $\alpha$-homogeneous
(\ref{krull}\,i).
Next, by \cite[Lemma 1.2(i)]{LP},
\begin{equation} \label{induced}
Q^G=\bigcap_{g\in G} Q^gkG \ ,
\end{equation}
and this intersection is finite. Thus, as right $kG$-module, $kG/Q^G$
is a finite subdirect product of the modules $kG/Q^gkG$. Furthermore,
each $kG/Q^gkG$ is generated, as $kG$-module, by the $k[H^g]$-submodule
$k[H^g]/Q^g$. Now $k[H^g]/Q^g$ is isomorphic to $kH/Q$ (as ring)
and is therefore $\alpha$-homogeneous. Consequently, each $kG/Q^gkG$ is
$\alpha$-homogeneous, by (\ref{krull}\,v), and hence so is $kG/Q^G$,
by (\ref{krull}\,i).

Finally, using \eqref{induced}, $I=\bigcap_{g\in G}(Q\cap kN)^g$ and,
by \cite[Lemma 4.2(i)]{B}, $kN/Q\cap kN$ is $\alpha$-homogeneous.
Therefore, so are its $G$-conjugates $kN/(Q\cap kN)^g$ and
consequently $kN/I$ as well, by the usual (finite) subdirect product
argument (\ref{krull}\,i). Since $kN/I$ generates the $kG$-module
$kG/IkG$, the latter is $\alpha$-homogeneous, by (\ref{krull}\,v), and
the proof is complete.
\end{proof}

\begin{prop}
Let $G$, $H$, $N$, and $Q$ be as in the above Lemma. Further, let $P$
be a prime ideal of $kG$ with $P\supset Q^G$. Then the following are
equivalent:
\begin{enumerate}
\item[i.] $P$ is minimal over $Q^G$;
\item[ii.] $P\cap kN=Q^G\cap kN$;
\item[iii.] $\hgt P = \hgt Q$.
\end{enumerate}
\end{prop}

\begin{proof}
(i $\Rightarrow$ ii) Put $I=Q^G\cap kN=\bigcap_{g\in G}(Q\cap kN)^g$,
as in the Lemma (and its proof), and observe that $I$ is a $G$-prime
ideal of $kN$.  Indeed, as in \eqref{cutdown}, $Q\cap kN=\bigcap_{h\in
H}Q_0^n$ for some prime $Q_0$ of $kN$, and so $I=\bigcap_{g\in
G}(Q\cap kN)^g=\bigcap_{g\in G}Q_0^g$.  Clearly, $I\subset P\cap
kN$. Suppose, by way of contradiction, that this inclusion is
strict. Then, by \cite[Lemma 3.4(i)]{LP3}, there exists an element
$c\in P\cap kN$ such that $c$ is regular modulo $I$. Since $kG/IkG$ is
free as left and right $kN/I$-module, $c$ is also regular in $kG$
modulo $IkG$.  Therefore, $\kdim(kG/IkG+ckG)_{kG} < \kdim(kG/IkG)$;
see \cite[13.6]{GW}.  Further, by the above Lemma,
$\kdim(kG/IkG)=\kdim(kG/Q^G)$. Since $IkG+ckG\subset Q^G+ckG$, we
obtain altogether that
$$
\kdim(kG/Q^G+ckG)_{kG} < \kdim(kG/Q^G)\ .
$$
By the above Lemma again, we know that $kG/Q^G$ is homogeneous. So the
last
inequality forces $c$ to be regular modulo $Q^G$ in $kG$; see
(\ref{krull}\,ii).
Inasmuch as $c\in P$, we deduce from \cite[10.8]{GW} that $P$ cannot be
minimal over $Q^G$, contrary to our hypothesis.\smallskip

(ii $\Rightarrow$ iii) Recall from the previous part of the proof that
$Q\cap kN=\bigcap_{h\in H}Q_0^n$ for some prime $Q_0$ of $kN$, and
$Q^G\cap kN=\bigcap_{g\in G}Q_0^g$. Therefore, (ii) implies that
$P\cap kN=\bigcap_{g\in G}Q_0^g$. By \cite[Lemma 2.5(i)]{LP}, we
conclude that $\hgt Q = \hgt Q_0 = \hgt P$, thereby proving (iii).
\smallskip

(iii $\Rightarrow$ i) Assume $\hgt Q = \hgt P$. If $P$ is not minimal
over $Q^G$ then $P$ contains a prime, $P'$, that is minimal over
$Q^G$.  Since we already have seen that (i) implies (iii), we know
that $\hgt P' = \hgt Q$ and, consequently, $\hgt P = \hgt P'$. But
this last equality is impossible, as $P$ properly contains $P'$. The
proof of the proposition is now complete.
\end{proof}

%%%%%%%%%%%%%%%%%%%%%%%%%%%%%%%%%%%%

\section{The proof} \label{proof}

\textit{In this section, $G$ will always denote a polycyclic-by-finite
group.
Thus, the group algebra $kG$ is
left and right noetherian and satisfies the strong second layer
condition.}

%=============================================
\subsection{The Roseblade subgroup} \label{rose}
Following Roseblade \cite{R}, a subgroup $H$ of $G$ is called
\textit{orbital} if $H$ has only finitely many $G$-conjugates. An
orbital subgroup $H$ is \textit{isolated} if there is no orbital
subgroup $H'$ of $G$ with $H\subsetneq H'$ and $[H':H]$
finite. Finally, a polycyclic-by-finite group is said to be
\textit{orbitally sound} if all its isolated orbital subgroups are
normal.  As in \cite{R}, we shall write
$$
\nio(G)
$$
for the intersection of the normalizers of all isolated orbital
subgroups of $G$.
The crucial properties of $\nio(G)$ for our purposes are the following:
\begin{itemize}
\item $\nio(G)$ is an orbitally sound normal subgroup of finite index
in $G$.  Moreover, $\nio(G)$ is the unique maximal such subgroup and
$\nio(G)$ contains every finite-by-nilpotent normal subgroup of $G$
(see \cite[1.3]{R}).
\item The group algebra of $\nio(G)$ is catenary \cite[Corrigenda]{R}.
\end{itemize}
Since images of orbitally sound groups are clearly orbitally sound as
well, the
maximality property of $\nio(G)$ stated above implies in particular
that, for every
epimorphism $\varphi: G\twoheadrightarrow\overline{G}$, one has
$\varphi(\nio(G))\subset
\nio(\overline{G})$.

The proof of the following lemma assumes some familiarity with
\cite{LP}, \cite{L}.

\begin{lem}
Assume that $G\neq\nio(G)$. Then every faithful prime ideal of $kG$
is induced from some proper subgroup of $G$ containing $\nio(G)$.
\end{lem}

\begin{proof}
Let $P$ be the prime in question, and let $H$ denote the normalizer of
its
vertex. Then $H$ contains $\nio(G)$ and $P$ is induced from $kH$; see
\cite[(2.2) and Lemma 2.5]{L}. If $H\neq G$ we are done; so assume that
$H=G$. In this case, \cite[Lemma 2.3]{L} implies that $P$ is induced
from
$k\Delta$, where $\Delta=\Delta(G)$ denotes the FC-center of $G$. Since
$\Delta$ is finite-by-abelian, $\Delta$ is contained in $\nio(G)$.
Therefore, by transitivity of induction \cite[Lemma 1.2(iii)]{LP}, $P$
is induced from $k[\nio(G)]$ as well and the proof is complete.
\end{proof}

%=============================================
\subsection{A special case} \label{special}
The following Lemma proves the requisite condition \eqref{E:ht} in the
special case where the lower prime is induced from a subgroup of finite
index
in $G$ whose group algebra is catenary.

\begin{lem} Let $H$ be a subgroup of finite index in $G$ such that $kH$
is catenary. Then, for any pair $P\subsetneq P'$ of adjacent primes of
$kG$ such that $P$ is induced from $kH$, one has $\hgt P' =\hgt P +
1$.
\end{lem}

\begin{proof} Write $P=Q^G$ for some ideal $Q$ of $kH$ which we may
assume prime,
by Lemma \ref{min}. By Proposition \ref{grouprings}, $\hgt Q = \hgt P$.
Choose an ideal $Q'$ of $kH$ containing $Q$ and maximal with
respect to $(Q')^G\subset P'$. Then, by (\ref{min}),
$Q'$ is prime and $P'$ is minimal over $(Q')^G$. Moreover,
Proposition \ref{grouprings} gives $\hgt Q' = \hgt P'$. Therefore,
it suffices to show that $\hgt Q' =\hgt Q + 1$ or, since
$kH$ is catenary, that $Q$ and $Q'$ are adjacent.
Suppose otherwise and fix a prime $I$ of $kH$ properly between $Q$ and
$Q'$.
Then $\hgt Q < \hgt I < \hgt Q'$. Further,
$P=Q^G\subset I^G\subset (Q')^G\subset P'$,
and so $P'$ contains a prime, $J$, of $kG$ minimal over $I^G$. By
Proposition \ref{grouprings}
again, $\hgt J = \hgt I$. Consequently, all inclusions in
$P=Q^G\subset J\subset P'$ are strict, contradicting adjacency of
$P$ and $P'$.
\end{proof}

%=============================================
\subsection{Proof of the Theorem} \label{pf}
To show that $kG$ is catenary, we argue by induction on the index
$[G:\nio(G)]$. The case $G=\nio(G)$ is known \cite{R}. So assume
$G\neq\nio(G)$, and let $P\subsetneq P'$ be adjacent primes of $kG$. We
have to show that $\hgt P' =\hgt P + 1$. In view of Lemma
\ref{heights}, we may pass to $G/P^\dag$ and thus reduce to the case
where $P$ is faithful. Note also that $[G:\nio(G)]\ge
[G/P^\dag:\nio(G/P^\dag)]$; see (\ref{rose}). Thus, by Lemma
\ref{rose}, $P$ is induced from $kH$ for some proper subgroup $H$ of
$G$ containing $\nio(G)$.  Since $\nio(G)$ is contained in $\nio(H)$,
by the maximality property of $\nio$, we have
$[H:\nio(H)]<[G:\nio(G)]$ and so $kH$ is catenary, by induction.
Lemma \ref{special} now completes the proof. \hfill\qed

%%%%%%%%%%%%%%%%%%%%%%%%%%%%%%%%%%%%

\section{Concluding remarks} \label{remarks}
\textit{$G$ will continue to denote a polycyclic-by-finite group.}
\medskip

Roseblade's work \cite{R} establishes catenarity of $kG$ up to
a finite index: $G_0=\nio(G)$ is a normal subgroup of finite index
in $G$ such
that $kG_0$ is catenary. Putting $R=kG_0$ and $S=kG$, we have
$$
S=R\*\G\ ,
$$
a crossed product of the \textit{finite} group $\G=G/G_0$ over $R$. It
is
tempting to approach the catenary problem for $kG$ from a more general
angle by establishing a positive answer to the following
\bigskip\newline
\textbf{Question.}
\textit{Suppose that $S=R\*\G$ is a crossed product with $R$ an affine
noetherian $k$-algebra
and $\G$ a finite group. If $R$ satisfies \eqref{E:ht}, will $S$ do so
as well?}
\bigskip\newline
Recall that (without any hypotheses on $R$) the ideals of $R$ of the
form $P\cap R$, for some prime ideal $P$ of $S=R\*\G$,
are precisely the $\G$-primes of $R$ and that these ideals are exactly
those of the form $Q_{\G}=\bigcap_{x\in\G} Q^x$ for some prime $Q$ of
$R$. Moreover, if $Q_{\G}\subset Q'_{\G}$ then we may assume, after
replacing $Q$ by a suitable $\G$-conjugate, that $Q\subset Q'$.
The latter will be adjacent if $Q_{\G}$ and $Q'_{\G}$ are adjacent
$\G$-primes. Finally, $P\cap R= Q_{\G}$ entails $\hgt P = \hgt Q$. See
\cite{LP2} for details. A positive answer to the above question would
therefore result if the following were true:
\bigskip\newline
\textbf{Question$^\prime$.}
\textit{Do adjacent primes of $S=R\*\G$ contract to adjacent $\G$-primes
of $R$?}
\bigskip\newline
Here is what appears to be known about Question$^\prime$:
\begin{itemize}
\item Without noetherian hypotheses, the answer is negative in
general.  Indeed, Kaplansky \cite{K} constructs an extension of
\textit{commutative} (non-noetherian) rings $R\subset S$, with
$S=R*\G$ and $\G$ the group of order 2, so that Question$^\prime$
fails for $S$.
\item Question$^\prime$ has a positive answer if $R$ is affine
commutative. This is an application of Schelter's catenarity result
for affine PI-algebras; see \cite[Proposition 3.3]{L}.
\end{itemize}
It would be interesting to have more definitive results about this
topic.

\begin{ack} The authors would like to thank K. A. Brown for
reading a preliminary draft of this article and for his
valuable comments.
\end{ack}

\end{document}